\documentclass[12pt, oneside,reqno]{amsart}   	% use "amsart" instead of "article" for AMSLaTeX format
\usepackage{geometry}
\usepackage{float}
\allowdisplaybreaks
\usepackage{graphicx}				% Use pdf, png, jpg, or eps with pdflatex; use eps in DVI mode
\usepackage{amssymb}
\usepackage{amsmath}
\usepackage[all]{xy}
\usepackage{mathtools}
\usepackage{tikz-cd}
\usepackage{enumitem}
\usepackage{tikz-cd}
\usepackage{appendix}
\usepackage[T1]{fontenc}
\usepackage{esint}

\usepackage[%
bookmarks=true,
colorlinks,
linkcolor=blue,
urlcolor=blue,
citecolor=blue,
plainpages=false,
pdfpagelabels,
final,
breaklinks=true\between
]{hyperref}
\usepackage{hyperref}
\usepackage[numbers,sort&compress]{natbib}

\newtheorem{thm}{Theorem}[section]
\newtheorem{lem}[thm]{Lemma}
\newtheorem{cor}[thm]{Corollary}
\newtheorem{prop}[thm]{Proposition}

\newtheorem{rem}{Remark}[section]

\numberwithin{equation}{section}

\newcommand{\A}{\mathcal{A}}

\def\Pb{\ifmmode{\Bbb P}\else{$\Bbb P$}\fi}
\def\Z{\ifmmode{\Bbb Z}\else{$\Bbb Z$}\fi}
\def\C{\ifmmode{\Bbb C}\else{$\Bbb C$}\fi}
\def\R{\ifmmode{\Bbb R}\else{$\Bbb R$}\fi}
\def\S{\ifmmode{S^2}\else{$S^2$}\fi}

\def\S{\cal S}

\def\HH{\mathcal{H}}

\newenvironment{pf}{\paragraph{Proof:}}{\hfill$\square$ \newline}

\begin{document}
\title[Knotted expanders]{On knotted self expanders of the mean curvature flow in $\R^3$}
\author{Alexander Mramor and Shengwen Wang}
\date{}

\address{Department of Mathematics, University of Oklahoma, Norman, OK 73019, USA}
\address{School of Mathematical Sciences, Queen Mary University of London, Mile End Road,
London E1 4NS}
\email{amramor@ou.edu, shengwen.wang@qmul.ac.uk}

\begin{abstract} In this article we produce many examples of knotted self expanders. For example, there is an expander isotopic to a trefoil knot with three half cylinders attached. On the other hand, we also show that sufficiently simple expanders, in particular asymptotically conical ones with one end or those homeomorphic to annuli, are unknotted. \end{abstract}
\maketitle

\section{Introduction}

In the mean curvature flow (abbreviated MCF), two of the most important categories of solitons are:

\begin{itemize}
    \item \textbf{Self shrinkers}, that is surfaces $\Sigma$ for which $\sqrt{-t} \Sigma$, $t \in (-\infty, 0)$ is a mean curvature flow; as $t$ increases $\sqrt{-t} \Sigma$ ``shrinks.'' These arise in the singularity analysis of the MCF when one parabolically dilates about a fixed point and capture much useful information about a forming singularity. 
   $\medskip$
    
    \item \textbf{Expanders}, that is surfaces $\Sigma$ for which $\sqrt{t} \Sigma$, $t \in (0, \infty)$ is a mean curvature; as $t$ increases $\sqrt{t}\Sigma$ ``expands.'' These flows model the expected evolution locally of a (weakly defined) mean curvature flow out of an asymptotically conical singularity i.e. for times just after it occurs. 
\end{itemize}

Another interesting and useful perspective on these solitons is that they are minimal surfaces in the conformal metrics $g^S_{ij} = e^{-|X|^2/2n} \delta_{ij}$ and $g^E_{ij} = e^{|X|^2/2n} \delta_{ij}$ respectively. By standard notions of curvature these metrics are ill-behaved, but it turns out they have signed \textbf{Bakry--Emery Ricci curvature}, where for a metric $g_{ij}$ ($\delta_{ij}$ above) and conformal factor $e^{-f}$ we have: 
\[ \text{Ric}_f =\text{Ric} + \text{Hess}(f)  \]
And so in this sense the shrinker metric, $g^S$, has positive weighted Ricci curvature and the expander metric, $g^E$, has negative weighted Ricci curvature. A natural question of course is to ask to what extent such weighted curvatures reflect properties of true Ricci positive or negative metrics and most optimistically one might hope that Bakry--Emery Ricci positive metrics share many of the same properties of the round metric on the sphere and Bakry--Emery Ricci negative metrics are akin to hyperbolic space. A natural avenue to investigate these questions is through their minimal surfaces.

A classical theorem of Lawson \cite{Lawson1970} says that minimal surfaces in round $S^3$ (or $S^3$ endowed with any other Ricci positive metric) must be unknotted, or in other words that for a fixed genus $g$ they must all be isotopic to the natural connect sum of $g$ tubular neighborhoods of round circles. This was confirmed for self shrinkers in the closed case by the authors in \cite{MramorWang2017} and was later extended to the noncompact case, see \cite{Mramor2020noncompact, ChuFranz2025}. On the other hand, there are minimal surfaces in hyperbolic space which are knotted, as constructed by Oliveira and Soret in \cite{deOliveiraSoret1998} (for a related paper see \cite{Martin2014Properly}); the point of this paper is to show that in line with the intuition above there are indeed knotted expanders in $\R^3$:
\begin{thm}\label{mainthm}
    For each $g \geq 1$ there are knotted expanders of genus $g$ with $r$ ends for $r$ large enough. For example, there is an expander isotopic to the boundary of a tubular neighbourhood of a trefoil knot, with three disks removed and a half-cylinder attached along each resulting boundary circle.
\end{thm}
 In fact we show obtain a more general existence statement than above, see theorem \ref{OS_analogy} below; other expanders with nontrivial genus were quite recently produced in \cite{ShaoZou2025Expanders, HoffmanMartinWhite2026}. By knotted in this generality we mean a surface which in large balls is not isotopic to a standardly embedded genus $g$ surface with a number of half-cylinders attached where some of them are potentially in nested families, although in the examples we produce the ends will all bound solid half-cylinders.

 An important quality of hyperbolic space for Oliveira and Soret is that there are minimal annuli/catenoids are relatively abundant in hyperbolic space and can be arranged in many different ways; it is perhaps less obvious that many such minimal annuli exist in the expander metric but sufficiently many expanding annuli for our purposes can be constructed via the Plateau problem using two correctly picked expanding discs attached by an appropriately sized cylindrical segment as area competitors. These can then used as barriers to show that appropriately constructed surfaces are incompressible in cleverly designed weakly mean convex domains so that there is no loss of topology when one minimizes in their isotopy class.
$\medskip$

On the other hand Olivera and Sorret can show some unknottedness results when the topology of a minimal surface in $\mathbb{H}^3$ is sufficiently simple (see corollaries 1,2 in \cite{deOliveiraSoret1998}), and indeed by analogy we can show closely related results in the expander case:

\begin{thm}\label{Unknottedness_Annuli}
Every connected, properly embedded, smoothly asymptotically conical self-expander in $\mathbb R^3$ that is homeomorphic to an annulus or has one end must be unknotted. 
\end{thm}
$\medskip$

To prove this result we borrow ideas from Meeks in \cite{Meeks1981}, showing that properly embedded, asymptotically conical expanders homeomorphic to a compact surface with finitely many punctures must be Heegaard splittings restricted to large balls about the origin. With this in hand the claim follows using Dehn's lemma and Waldhausen's theorem respectively. 

\textbf{Acknowledgements:} This project begun in August 2024 during the first named author's visit to QMUL and he thanks the department for their hospitality. 
$\medskip$

\textbf{Statement on AI use:} AI was helpful in confirming our (human generated) ansatz for the expanding annuli constructed below, both in the area estimate and in pointing out the important reference \cite{HoffmanMartinWhite2026}. It was also used in proofreading. 

\section{Construction of expanding annuli as barriers}

\begin{figure}[h]
    \centering
    \includegraphics[width=0.65\linewidth]{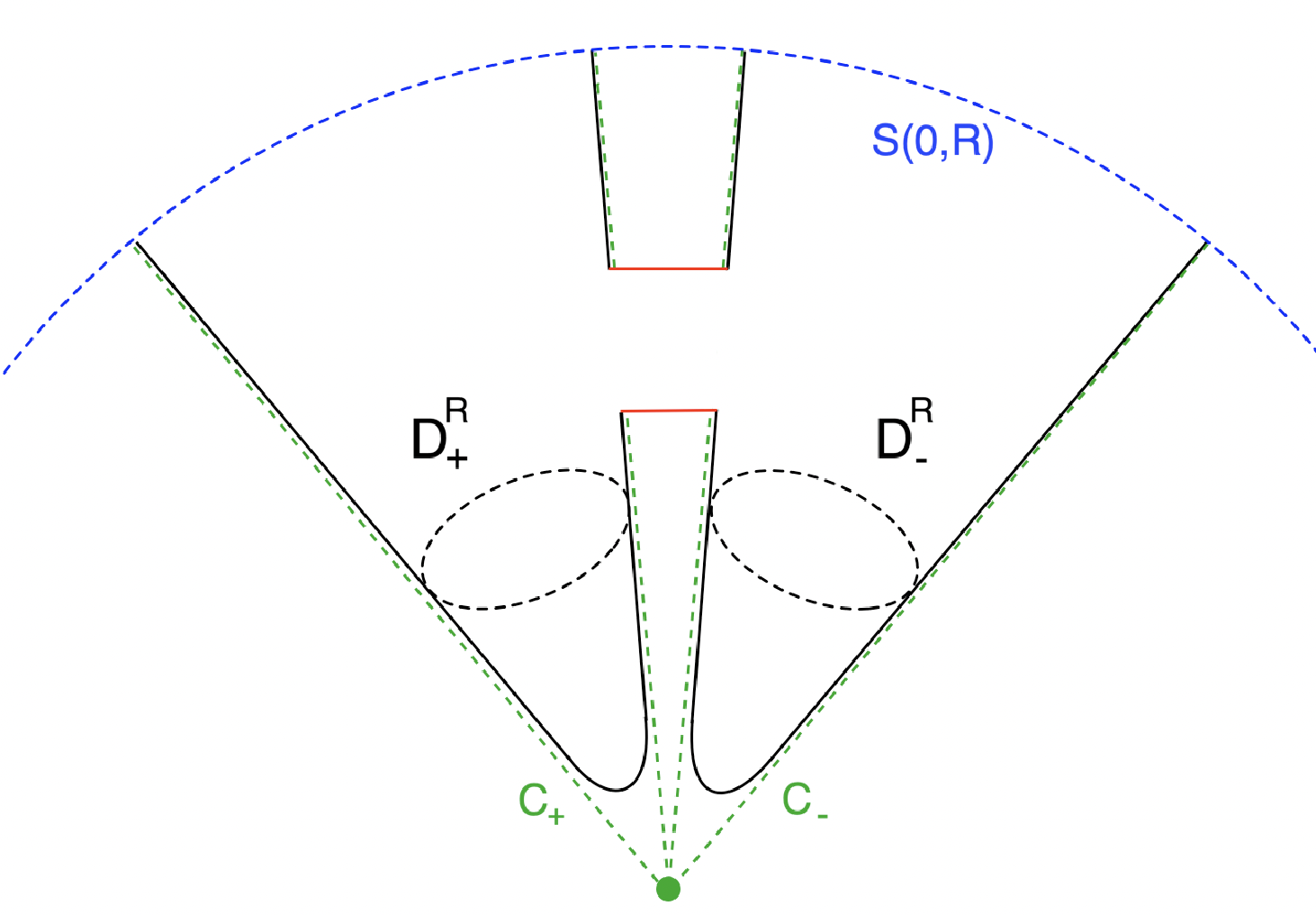}
    \caption{An illustration of the bridged together expanding discs $D_\pm$ considered in a large ball $B(0,R)$ for use as a competitor in the Plateau problem. The bridge, indicated in red above, is given by a short rectangular box; for illustrative purposes its relative size is exaggerated and also, as we'll see below, the sides of $C_\pm$ facing each other will also typically need to be ``flattened.'' }
    \label{cutpaste}
\end{figure}

We recall that the weighted area functional $$\A(S)=\int_S e^{|X|^2/4}\,d\HH^2.$$
is the area in the expander metric
$g=e^{|X|^2/4}g_{\mathrm{Euc}}$; minimal surfaces, satisfying $H + \frac{\langle X, \nu \rangle}{2} = 0$, are called expanders. As the title suggests, our goal in this section is to construct annuli which are minimal surfaces in this weighted area. To construct these we will: 
\begin{enumerate}
    \item Construct expanding annuli $A_R$ in $B(0,R)$ by the direct method applied to the area functional. To produce a competitor we will take two expanding discs $D_\pm$ and attach them via a bridge as in figure \ref{cutpaste}. As discussed more below simple reasoning suggests that one should be able to arrange such a bridge to decrease expander area, at least in many cases, which we'll proceed to check carefully. Consequently, we obtain an expanding annulus $A_R$ in $B(0,R)$. 
    $\medskip$
    
    \item Take a limit of the surfaces $A_R$ as $R \to \infty$. Since each of the $A_R$ are area minimizing they are stable so we may apply classical curvature estimates, and we can also produce area bounds. To ensure the limit is nontrivial (or roughly speaking that the joining neck doesn't degenerate somehow) we will apply a recent observation of Hoffman--Martin--White \cite{HoffmanMartinWhite2026} which roughly says that critical points of the distance functional along the smooth expanders we produce must be a uniformly bounded distance from the origin.
\end{enumerate}

We first discuss in more detail the discs $D_\pm$ we use, at least preliminarily. For the sake of simplicity in the calculations below we start with choosing a round cone $C$ of small opening angle and subsequently fixing an $R_0$ and some $R \gg R_0$ (so one can apply necessary estimates in the following). With these fixed, we next consider two copies of $C$ and rotate one of them about the $y$-axis so that their links are disjoint and apart by an angle over $R_0^{-3}/2$. Denoting the resulting cones by $C_\pm$, by \cite{Ding2020Minimal} there are 2 inward expander disks $D^R_+,D^R_- \subset B(0,R)$ of least weighted area that agree with $C_\pm$ along $ S(0,R)$. Note that they will be disjoint because the cones $C_\pm$ bound disjoint domains. 

\begin{rem}[Ansatz]
Because $C_\pm$ and hence $D_\pm$ are disjoint but very close intuitively cutting out relatively large (compared to the separation) discs from either of them and attaching a short bridge between them should decrease area. The situation turns out to actually be a bit more complicated than this essentially due to curvature of the link, but we confirm this in sufficiently many cases for our applications below. 
\end{rem}

By lemma 5.9 of \cite{Ding2020Minimal} and that the discs $D^R_\pm$ converge as $R \to \infty$, we see that for $R \gg R_0$ at distance $R_0\leq r \leq 2R_0$ from the origin the disks $D^R_\pm$ are graphical over their asymptotic cone with norms of order $\frac{k_g}{r}+O(\frac{1}{r^3})$ and gradient or order $O(\frac{1}{r^2})$, where $k_g$ is the geodesic curvature of link. Up to rotation we may suppose that the axis of the original cone $C$ is the $z$-axis and that in defining $C_\pm$ we have that over the rectangle $Q_{R_0}=[-50/R_0,50/R_0] \times [R_0-10/R_0,R_0]$ in the $yz$-plane the facing patches of $D^R_\pm$ (independent of $R$!) are graphs $x=u_\pm(y,z)$ satisfying, uniformly on $Q_{R_0}$,

\[ u_\pm=\pm \frac{k_g}{R_0}+O(\frac{1}{R_0^2}), \text{ }|Du_\pm|=O(R_0^{-2}) \]
where in the above the rotational directions of $C_\pm$ gave us error terms of order $R_0\cdot\frac{1}{R_0^3}=\frac{1}{R_0^2}$.

The patches are symmetric by reflection and so their boundaries can be connected by straight lines in the following gluing procedure: remove these two graph patches, denoted by $P_\pm$, and join their boundaries by the four-sided wall $$ B_{R_0}=\{(x,y,z):(y,z)\in\partial Q_{R_0},\   u_-(y,z \leq x\leq u_+(y,z)\}.$$
Denote the resulting surface by $N_R$. The size of the cubes $Q_{R_0}$ implies $$|y|,|z-R_0|=O(\frac{1}{R_0})$$ and the small graphicality bound of patches implies $$|x|=O(\frac{1}{R_0}).$$ 
Thus:
\begin{align*}
 \frac{|X|^2}{4}=&\frac{1}{4}[x^2+y^2+z^2]\\
 =&\frac{1}{4}[x^2+y^2+(z-R_0)^2+2(z-R_0)R_0+R_0^2]\\
 =&\frac{x^2+y^2+(z-R_0)^2}{4}+\frac{R_0^2}{4}+\frac{R_0(z-R_0)}{2}\\
 =&\frac{R_0^2}{4}+\frac{R_0(z-R_0)}{2}+O(R_0^{-2}).
\end{align*}
We also notice that:
\begin{align*}
 e^{|X|^2/4}=&e^{R_0^2/4}e^{R_0(z-R_0)/2}e^{O(R_0^{-2})}\\
 =&e^{R_0^2/4}e^{R_0(z-R_0)/2}\left(1+O(R_0^{-2})\right)\\
 \sqrt{1+|Du_\pm|^2}=&1+o(1).
\end{align*}
As a consequence, the two removed patches satisfy: 
\begin{align*}
 \A(P_+)+\A(P_-) &=\sum_{\pm}\int_{Q_{R_0}}e^{(u_\pm^2+y^2+z^2)/4} \sqrt{1+|Du_\pm|^2}\,dy\,dz\\
 &=e^{R_0^2/4}(2+o(1))\frac{100}{R_0}  \left[\frac{2}{R_0}e^{R_0(z-R_0)/2}\right]_{R_0-10/R_0}^{R_0}\\
 &=e^{R_0^2/4}(2+o(1))\frac{100}{R_0}  \frac{2}{R_0}(1-e^{-5})\\
 &=\frac{e^{R_0^2/4}}{R_0^2}   \left[400(1-e^{-5})+o(1)\right].
\end{align*}
Here $100/R_0$ is the length of the $y$-interval of $Q_{R_0}$.

Each wall face lies in a coordinate plane, so its Euclidean area element is $dx\,d\ell$, where $d\ell$ is arclength along the corresponding side of $\partial Q_{R_0}$. Since $u_+-u_-=2{k_g}/R_0+o(R_0^{-1})$, we have
\begin{align*}
 \A(B_{R_0}) &=\int_{\partial Q_{R_0}}\int_{u_-}^{u_+} e^{(x^2+y^2+z^2)/4}\,dx\,d\ell\\
 &=e^{R_0^2/4}\left(\frac{2{k_g}}{R_0}+o(R_0^{-1})\right) \int_{\partial Q_{R_0}}e^{R_0(z-R_0)/2}\,d\ell\\
 &=e^{R_0^2/4}\left(\frac{2{k_g}}{R_0}+o(R_0^{-1})\right) \left[ \int_{-50/R_0}^{50/R_0}(1+e^{-5})\,dy +2\int_{R_0-10/R_0}^{R_0}e^{R_0(z-R_0)/2}\,dz \right]\\
 &=e^{R_0^2/4}\left(\frac{2{k_g}}{R_0}+o(R_0^{-1})\right) \left[ \frac{100}{R_0}(1+e^{-5}) +2\left[\frac{2}{R_0}e^{R_0(z-R_0)/2}\right]_{R_0-10/R_0}^{R_0}  \right]\\
 &=e^{R_0^2/4}\left(\frac{2{k_g}}{R_0}+o(R_0^{-1})\right) \left[ \frac{100}{R_0}(1+e^{-5}) +\frac{4}{R_0}(1-e^{-5}) \right]\\
 &=\frac{e^{R_0^2/4}}{R_0^2}{k_g} \left[200(1+e^{-5})+8(1-e^{-5})+o(1)\right].
\end{align*}
The first term in the square brackets comes from the two sides $z=R_0$ and $z=R_0-10/R_0$, each of length $100/R_0$; the second comes from the two sides $y=\pm50/R_0$.

The difference between the two limiting coefficients is $$400(1-e^{-5})-k_g\left[200(1+e^{-5})+8(1-e^{-5})\right]$$ 
We see that if $k_g$ on $C$ was sufficiently small (less than one works, for example) this would be positive so that we would have the connect sum above would decrease area, which is what we want, but this is restriction is less than optimal. It might be possible to instead use slightly different parameters in the construction above, but to sidestep this issue we instead deform the original cones $C_+$ and $C_-$ so that along the sides facing each other the curvature of their links is small, and in fact the link is given locally by a piece of a round cone of small curvature. If this ``flattening" of $C_+$ and $C_-$ is done in a smooth, convex way (which can certainly be arranged) on each the area minimizing expanders asymptotic to them (of course, considered separately) from \cite{Ding2020Minimal} will be mean convex, and since they are expanding this will imply they are discs $D_{\pm}$ (although if they were topologically nontrivial for our purposes it wouldn't be so important). By the convergence, this also implies the weighted area minimizing solutions $D_{\pm}^R$ are also topologically discs. As before, they will continue to be disjoint because the cones bound disjoint domains. 
$\medskip$

As an upshot, for such modified $C_+, C_-$ the computation above implies for sufficiently large $R$ and every
ball $B_s$ containing the surgery region, $$ \A(N_R\cap B_s) <\A(D^R_+\cap B_s)+\A(D^R_-\cap B_s).$$

This is analogous to the catenoid estimate in \cite{KetoverMarquesNeves2020}. The corners and edges can be smoothed with arbitrarily small change in weighted area, preserving this strict inequality. Recalling that $R \gg R_0$ and that the discs $D^R_\pm$ can be taken to be of least area then by minimizing in the isotopy class of $N_R$ via Meeks and Yau \cite{Meeks1982TheEO} we obtain the following: 
\begin{lem}
     For $R$ sufficiently large, there exists an expanding annulus $A_R \subset B(0,R)$ which locally minimizes the expander area. 
\end{lem}
We then wish to take a (subsequential) limit of these via Arzela-Ascoli along a sequence $R_i \to \infty$. By the geometric Arzela-Ascoli theorem and regularity theory for minimal surfaces we need to check:
\begin{enumerate}[label=\alph*)]
    \item We have sufficient curvature bounds on compacta -- this is standard because each of the $A_{R_i}$ are area minimizing.
    \item The areas of the $A_{R_i}$ are uniformly bounded on a given bounded domain.
    \item The limit, provided it exists from items (a) and (b), is an annulus itself. 
\end{enumerate}
An important ingredient for us is the following recent result of Hoffman, Martin, and White: 
\begin{thm}[Theorem 12.4 in \cite{HoffmanMartinWhite2026}]\label{HMW-Thm}
Suppose that $C$ is a cone that is smooth except at the origin, and
suppose that $0 < \varepsilon < \pi/2$. There exists $\tilde{R} < \infty$,
depending only on $C$ and $\varepsilon$, with the following property.
If $R > \tilde{R}$ and if $E$ is an expander in $B(0,R)$ such that
\[
  \partial E = C \cap \partial B(0,R),
\]
then
\begin{enumerate}
  \item $E \setminus B(0,\tilde{R})$ is in the interior of the domain of
        the nearest point retraction $\pi(\cdot)$ onto $C$.
  \item $\pi$ maps $E \setminus B(0,\tilde{R})$ diffeomorphically onto a
        subset of $C$.
  \item For each $p \in E \setminus B(0,\tilde{R})$, the angle between
        $\operatorname{Tan}(E,p)$ and $\operatorname{Tan}(C,\pi(p))$ is
        $<\varepsilon$, and the angle between $\operatorname{Tan}(E,p)$
        and $p$ is $<\varepsilon$.
\end{enumerate}
\end{thm}
Since the expanders $A_{R_i}$ arose as solutions to the Plateau problem with boundary given by $(C_+ \cup C_-) \cap S(0,R_i)$ and $C_+ \cup C_-$ is a smooth cone away from the origin, we can apply parts (1) and (2) of the result above to see that outside some large uniform ball $B(0, \tilde{R})$ each of the $A_{R_i}$ project onto $(C_+ \cup C_-)$. By cutting and pasting techniques (say, using a fixed $N_R$ and deforming it slightly along the sphere $r = \tilde{R}$) its easy to produce area competitors of uniformly bounded area in $B(0, \tilde{R})$, so combined with the above each of the $A_{R_i}$ have uniformly bounded area on bounded domains. Part (3) of the theorem above also implies that the ``bridges'' of the annuli $A_{R_i}$ cannot degenerate to spatial infinity, because these have corresponding points where the tangent plane is perpendicular to the radial vector field.

We recall that the links of the cones $C_{\pm}$ we started with were originally round, but we could it could be the case they would have to be "flattened" on a side somewhat. So to summarize we have the following, stated somewhat imprecisely but sufficient for our purposes: 
\begin{prop}\label{expanding_cat} Suppose that $c_1, c_2$ are two disjoint round circles on $S^2$ of the same radius. If they are sufficiently close together there is an asymptotically conical expander $A_{c_1, c_2}$, homeomorphic to an open annulus, whose link is $C^0$ close to $c_1 \cup c_2$.     
\end{prop}

%In fact, and important for the sequel, it’s not hard to see we can arrange disjoint families of expanding annuli. If we consider a collection of $k$ pairs of disjoint cones which bound pairwise disjoint domains with each pair fulfilling the above, the expanding discs asymptotic to them will be pairwise disjoint as well, and its easy to see we may produce pairwise disjoint annular competitors reasoning as before. The argument above then gives $k$ disjoint expanding annuli. 

The links of two small cones and the construction of annuli can be carried out inside an arbitrarily small round cone with vertex at origin. So the expanding annuli are inside an arbitrarily small round cone. Since one can get as many such round cones as one needs, we obtain $k$ pairwise disjoint expanding annuli for any finite $k$.

\section{Producing knotted expanders}

With the expanding annuli above in hand we may easily adapt the construction of de Oliveira and Soret \cite{deOliveiraSoret1998} to the expander metric, we obtain the following theorem, which implies Theorem \ref{mainthm}.
\begin{thm}[c.f. theorem 2 in \cite{deOliveiraSoret1998}]\label{OS_analogy}
Suppose that $\Sigma \subset \R^3$ is an embedded and orientable surface which bounds a handlebody $H_\Sigma$. Then for $r$ sufficiently large there exists a complete properly embedded self-expander $M$, stable for the weighted area $\A$, with $r$ ends and bounding a region $H_M$ such that the interior of $H_\Sigma$ is isotopic to the interior of $H_M$ in balls $B(0,R)$ for $R$ sufficiently large.    
\end{thm}

\begin{pf}
The expanding annuli from Proposition \ref{expanding_cat} replace the catenoids in the construction of de Oliveira and Soret to prevent compression. In the expander metric $g=e^{|X|^2/4}g_{\mathrm{Euc}}$, the annular barriers are minimal and the spheres $\partial B_R$ are strictly mean-convex. Choosing $R$ so that $\partial B(0,R)$ meets the annular barriers transversely, the truncated domain has piecewise smooth mean-convex boundary with interior dihedral angles less than $\pi$. Thus the Meeks--Yau minimization applies in the compact domains obtained by intersecting the exterior of the fixed barriers with $\overline{B(0,R)}$, producing expanders $\Sigma_R$ in the prescribed isotopy class. Choosing their boundaries consistently as $C\cap\partial B(0,R)$ for a fixed smooth cone $C$, theorem \ref{HMW-Thm} give uniform graphical annular ends and area bounds as in the proof proposition \ref{expanding_cat} so we can
choose a sequence $R_i\to\infty$ along which $\Sigma_{R_i}$ converges smoothly and with multiplicity one to a properly
embedded self-expander. So, all the ingredients for the topological construction of de Oliveira and Soret are available here, and the same construction proves the general existence assertion of theorem \ref{OS_analogy}. \end{pf}

A more precise statement on the number of ends $r$ can be given in terms of the ``bridge number," see \cite{deOliveiraSoret1998}. Because the isotopy of $\Sigma$ is not specified, as a consequence we get many knotted examples of expanders. In the remainder of this section, we illustrate how the construction works for a specific example in a hands-on way: 
\begin{cor}
    There exists an expander which, in balls $B(0,R)$ for $R$ sufficiently large, is isotopic to the boundary of tubular a neighborhood of a trefoil knot with three cylindrical ends attached. 
\end{cor}

\begin{figure}[h]
    \centering
    \includegraphics[width=1\linewidth]{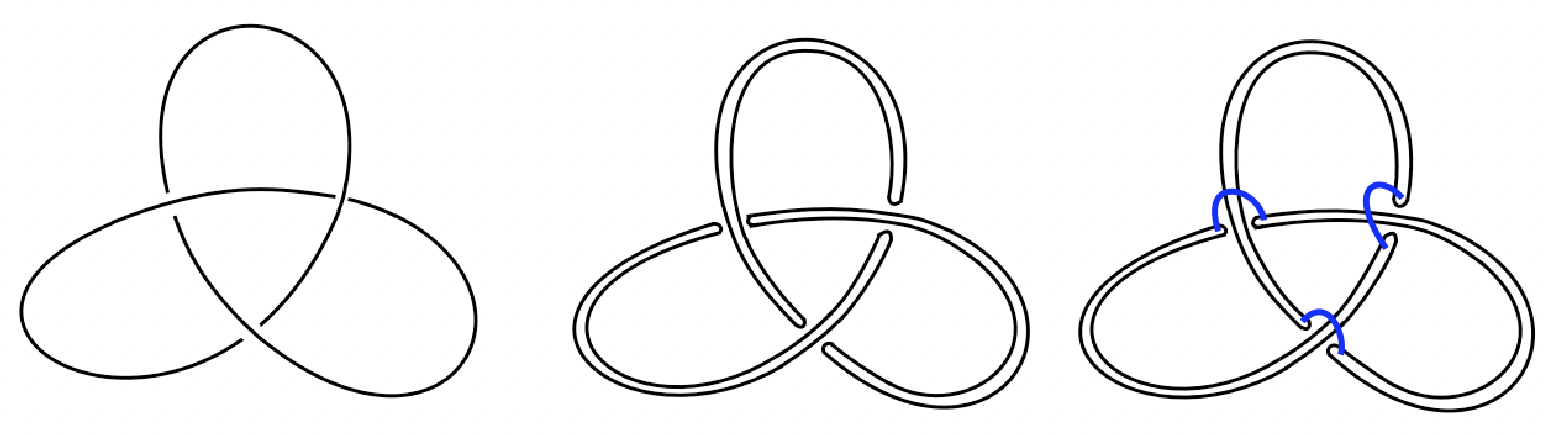}
    \caption{From left to right we have a projection of a trefoil knot $K$ onto $S(0,R)$ (where $R$ is so large in comparison to our representative for $K$ it might as be $\R^2$ for our purposes), the surgered thickening of the projection of $K$, and lastly the overpasses of $K$ reinstalled via bridges that are topologically cylindrical segments. Here we see that the final surface has three boundary components, corresponding to 3 ends in the limit.}
    \label{trefoil}
 
\end{figure}
\begin{pf}
    Consider a large ball $B(0,R)$ and let $\gamma$ be an embedding of a trefoil knot in it for which the projection $\tilde{\gamma}$ of $\gamma$ onto $S(0,R)$ is immersed (but fails to be embedded where it crosses over itself). Then we can "thicken" $\tilde{\gamma} \subset S(0,R)$, surgering along the thickened $\tilde{\gamma}$ depending on if a crossing is over or under, and finally add bridges facing into interior of $B(0,R)$ for over crossings as indicated in figure \ref{trefoil}. 

    Denoting this surface by $\Sigma_R$, an important issues that it is not incompressible; when we attempt to solve the Plateau problem, it could be the case that the bridges (blue in the figure) break. To remedy this, we use can use the expanding annuli from proposition \ref{expanding_cat} as inner barriers; by adjusting the thickenings of $\tilde{\gamma}$ and the size of the bridges appropriately, expanding annuli can be slipped inside each of them. Continuing to denote this surface by $\Sigma_R$, removing these expanding annuli from $B(0,R)$ to get an expander mean convex domain $U$ within which $\Sigma_U$ lays as an incompressible surface. For each $R$ sufficiently large then, we can produce an expander $\Sigma_R \subset B(0,R)$ which is isotopic to a trefoil knot with 3 half cylinders attached. Choosing the link of these boundary components consistently to coincide with smooth cones we can then proceed as the in proof of proposition \ref{expanding_cat} to obtain our example as a limit of $\Sigma_{R_i}$, $R_i \to \infty$. 
    \end{pf}

The general situation is much the same: if a surface $\Sigma \subset B(0,R)$ bounds a handlebody $H_\Sigma$, then we can retract it onto its graph $G$ and consider its projection $\tilde{G}$ onto $S(0,R)$. As in the example above, there might be some crossings where bridges are needed and as above we may use expanding annuli as barriers. However there can be extra complications: for example there might also be cycles of $\tilde{G}$ which don't come from self crossings, and so even after adding in overpasses there could be nested boundary components of the thickening of $\tilde{G}$. To see how this can be an issue, consider for instance the simple case $\tilde{G}$ is a round circle, corresponding to starting with a round torus. Then when one fattens $\tilde{G}$ the boundary is given by two concentric round circles; the solution to the plateau problem in $B(0,R)$ in this case is apt to be an annulus (or worse, a pair of discs), not a torus with some number of punctures which is what theorem \ref{OS_analogy} produces. To deal with this, we may arbitrarily add a catenoidal overpass along the thickening of $\tilde{G}$ in $B(0,R)$; the boundary would then be connected then and the starting surface is a once punctured torus. An issue is that this is compressible; to deal with this, one can use a catenoid linking with $\tilde{G}$ as indicated in figures 2.2.2 and 3.1.2 of \cite{deOliveiraSoret1998}.

\section{Unknottedness of topologically simple expanders}

In the following, a Heegaard splitting of a 3-ball $B$ is an embedded surface in $B$ which decomposes it into two handlebodies. 

\begin{lem}\label{No_Max}
Let $M\subset\mathbb R^3$ be a properly embedded, smoothly asymptotically conical surface which is homeomorphic to a compact connected surface punctured in finitely many points. If, for some $a\in\mathbb R^3$,
the function $|X-a|^2$ restricted to $M$ is Morse with no local maxima, then $M \cap B(a,R)$ for $R \gg 0$ is a Heegaard splitting of $B(a,R)$.
\end{lem}

\begin{proof}
By applying Morse theory to the distance squared $f_a = |X -a|^2$ to the point $a$, the argument for theorem 3 on page 406 of \cite{Meeks1981} gives the two components of $\R^3$ bounded by $M$ are diffeomorphic to regular neighborhoods of 1-dimensional proper CW-complexes. Denoting these domains $I$ and $O$, the idea is that because there are no index 2 critical points, only 0- and 1-handles are added to $I$ and $O$ through critical points of $f_a$ along $M$. Noting by the asymptotically conical assumption that $f_a$ must have only finitely many critical points, its clear that in bounded balls both $I$ and $O$ retract to finite graphs; because $M$ is connected and asymptotically conical for $R$ sufficiently large $I \cap B(a,R)$ and $O \cap B(a,R)$ will be connected. From this for $R$ large $I \cap B(a,R)$ and $O \cap B(a,R)$ will be handlebodies and $M \cap B(a,R)$ is a Heegaard splitting of $B(a,R)$ for $R \gg 0$ as claimed. 
\end{proof}

\begin{proof}[Proof of Theorem~\ref{Unknottedness_Annuli}] We recall that for generic $a$, $f_a=|X-a|^2$ Morse on $M$; choose a generic center $a$ arbitrarily close to the origin. 
Using $\Delta_M X=\frac12X^\perp$, we obtain
\begin{align*}
\Delta_M f_a =4+|X^\perp|^2-\langle a,X^\perp\rangle \geq 4-\frac{|a|^2}{4}>0
\end{align*}
when $|a|<4$; thus $f_a$ has no local maxima. With this in mind we then argue essentially as Meeks does as in his proof of the corollary on page 408 of \cite{Meeks1981}, where he shows an analogous statement for minimal surfaces in $\R^3$ endowed with a metric of negative sectional curvature as we next discuss.

Starting with the annulus case, we first note that $M$ is diffeomorphic to a twice punctured sphere so the lemma above can be applied. For $R \gg 0$ one of the domains $M \cap B(a,R)$ bounds is topologically a solid cylindrical segment $C$ with two boundary components $\sigma_1, \sigma_2$ on $S(a,R)$; denote the core of $C$ by $\gamma$ which we can push across $C$ to be in $\partial C$ with boundary points $p_1, p_2$ on $\partial\sigma_1, \partial\sigma_2$ respectively. The endpoints of $\gamma$ can then be joined together by a curve in $S(a,R) \setminus C$ to get a simple closed curve $\widetilde{\gamma}$ along the boundary of $\overline{B(a,R) \setminus C}$. Noting that by Mayer-Vietoris $H_1(\overline{B(a,R) \setminus C}) \simeq \Z$, and that because $M \cap B(a, R)$ is a Heegaard splitting $\overline{B(a,R) \setminus C}$ is a handlebody, we see it in fact must be a solid torus; from this its easy to see $\widetilde{\gamma}$ is nullhomotopic in it. Thus by Dehn's lemma it bounds a disc which implies that $M$ is unknotted.

 Considering next the case $M$ has one, asymptotically conical, end its easy to see its diffeomorphic to a one punctured compact surface $C$ so again the lemma above applies. Consequently we have that $M \cap B(a,R)$ is a Heegaard splitting of $B(a,R)$ with one boundary component for $R$ sufficiently large, so Waldhausen's theorem \cite{waldhausen1968heegaard} implies it must be unknotted. In the one ended case one can also proceed (as Meeks originally does) via ttheorem 1 in \cite{Meeks1981}, shown using Morse theory, as indicated in the proof of corollaries 2, 4 therein. 
\end{proof}

\begin{rem}
    Note that the argument above also implies that the knotted expanders we produced in the previous section are also Heegaard splittings in large balls, but as Hall's constructions \cite{Hall1984} have already shown this is generally not enough to imply that they are unknotted. There is a closely related but stronger notion of Heegaard splitting in terms of fundamental groups which does imply a surface is unknotted however even if it has many ends, as shown in Frohman \cite{frohman1992heegaard} and utilized in \cite{ChuFranz2025}.
\end{rem}

\bibliographystyle{plain}
\bibliography{bibliography}

\end{document}